\documentclass[12pt]{article}
\usepackage{amsfonts}
\usepackage{latexsym}

\newcommand{\be}{\begin{equation}} \newcommand{\ee}{\end{equation}}
\newcommand{\bea}{\begin{eqnarray}} \newcommand{\eea}{\end{eqnarray}}
\newcommand{\bean}{\begin{eqnarray*}}
  \newcommand{\eean}{\end{eqnarray*}}
\newcommand{\brray}{\begin{array}} \newcommand{\erray}{\end{array}}
\newcommand{\ben}{\begin{equation}{nonumber}}
  \newcommand{\een}{\end{equation}{nonumber}}
\newcommand{\newsection}[1]{\setcounter{equation}{0}
  \setcounter{dfn}{0}
\section{#1}}

\newtheorem{dfn}{Definition}[section] \newtheorem{thm}[dfn]{Theorem}
\newtheorem{lmma}[dfn]{Lemma} \newtheorem{ppsn}[dfn]{Proposition}
\newtheorem{crlre}[dfn]{Corollary} \newtheorem{xmpl}[dfn]{Example}
\newtheorem{rmrk}[dfn]{Remark}

\newcommand{\bdfn}{\begin{dfn}} \newcommand{\bthm}{\begin{thm}}
    \newcommand{\blmma}{\begin{lmma}}
      \newcommand{\bppsn}{\begin{ppsn}}
        \newcommand{\bcrlre}{\begin{crlre}}
          \newcommand{\bxmpl}{\begin{xmpl}}
            \newcommand{\brmrk}{\begin{rmrk}}

              \newcommand{\edfn}{\end{dfn}}
            \newcommand{\ethm}{\end{thm}}
          \newcommand{\elmma}{\end{lmma}}
        \newcommand{\eppsn}{\end{ppsn}}
      \newcommand{\ecrlre}{\end{crlre}}
    \newcommand{\exmpl}{\end{xmpl}} \newcommand{\ermrk}{\end{rmrk}}

\newcommand{\IC}{{\Bbb C}}

 \newcommand{\IN}{{\Bbb N}}
 
 \newcommand{\IR}{{\Bbb R}}
 \newcommand{\IT}{{\Bbb T}}

 \newcommand{\IZ}{{\Bbb Z}}

 \newcommand{\al}{\alpha}

\newcommand{\dlt}{\delta} 
 
\newcommand{\lmd}{\lambda}

 \newcommand{\del}{\partial}

\newcommand{\cla}{{\cal A}} \newcommand{\clb}{{\cal B}}
 \newcommand{\cld}{{\cal D}}
 
 \newcommand{\clh}{{\cal H}}
 
 \newcommand{\cll}{{\cal L}}

\newcommand{\AAZ}{{\cal A}_{\hbar{}}^{\infty}}

  \def 
\bbE
{\mbox{\boldmath $E$}}                     

  \def 
\bbe
{\mbox{\boldmath $e$}}

 \def\a*{{\cal A}_{h,*}} \def\B{{\cal B}(h)}
\def\B1{{\cal B}_1(h)} \def\b{{\cal B}^{s. a. }(h)} \def\b1{{\cal
    B}^{s. a. }_1(h)}

\newcommand{\raro}{\rightarrow}

  \def \qed 
{
  \mbox{}\hfill $\Box$\vspace{1ex}} 

\begin{document}

\author{{\large \sc Partha Sarathi Chakraborty}\\
         Indian Statistical
Institute,\\[-.5ex]
         203, B. T. Road,
Calcutta--700\,035, INDIA\\[-.5ex]
         email:
parthasc$\underline{\mbox{\hspace{.5em}}}$r@isical.ac.in\\[1ex]
{\large \sc Kalyan Bidhan Sinha}\\
         Indian Statistical
Institute,\\[-.5ex]
         7, SJSS Marg, New Delhi--110\,016,
INDIA\\[-.5ex]
         email: kbs@isid.ac.in}
\title{Geometry on the Quantum Heisenberg Manifold}
\maketitle

\section {Introduction}
Let G be the Heisenberg group. \\
$G =  \{  \left(  \matrix { 1 & x & z \cr 0 & 1 & y \cr 0 & 0 & 1 \cr } \right)  | x,y,z \in \IR \}$ \\
For a positive integer c, let $H_c$ be the subgroup of $G$ obtained when $ x,y,cz$ are integers. The Heisenberg manifold $M_c$ is the quotient $ G/H_c$. Nonzero Poisson brackets on $M_c$ invariant under left translation by $G$ are parametrized by two real parameters $ \mu,\nu$ with $ \mu^2 + \nu^2 \ne 0$ \cite {RI1}. For each positive integer $c$ and real numbers $\mu,\nu$ Rieffel constructed a C*-algebra $A_{\mu,\nu}^{c,\hbar}$ as example of deformation quantization along a Poisson bracket \cite {RI1}. These algebras have further been studied by \cite {AB1}  \cite {AB2} \cite {W1}. It was also remarked in \cite {RI1} that it should be possible to construct example of non-commutative geometry as expounded in \cite {CON}  in these algebras also. It is known \cite {RI1} that  Heisenberg group acts ergodically on $A_{\mu,\nu}^{c,\hbar}$  and $A_{\mu,\nu}^{c,\hbar}$ accomodates a unique invariant tracial state $\tau$. Using the group action we construct a family of spectral triples. It is shown that they induce same element in K-homology. We also show that the associated Kasparov module is non-trivial. This has been achieved by constructing explicitly the pairing with a unitary. We also compute the space of forms as described in \cite {CON} \cite {FGR}. Then we characterize torsionless and unitary connections. From that easily follows that a torsionless unitary connection can not exist. For a family of unitary connections we compute Ricci curvature and scalar curvature as introduced in \cite {FGR}. This family has non-trivial curvature. Following \cite {SAV} one can construct Quantum Dynamical Semigroups which will be natural candidate for heat semigroup. In a forthcoming version of the paper we wish to describe the construction and
dilate the semigroup.

Organization of the paper is as follows. In section 2 after introducing the algebra we compute the GNS space of  $\tau$ using a crucial result of Weaver \cite {W1}. In the next section following a general principle of construction of spectral triple on a C*-dynamical system with dynamics governed by a Lie group we construct spectral triples and compute the hyper trace  \cite {VAR} , \cite {CGS}  associated with the spectral triple.In section 3 we compute the space of forms (Chapter V of \cite {CON}). There are not too many instances of this computation in the literature.In section 4 after briefly recalling the notions introduced in \cite {FGR} we compute the space of $L^2$-forms. Then we characterize torsionless/unitary connections and show a connection can not simultanaeously be torsionless and unitary. In the next section for a concrete family of unitary connections we compute Rici curvature and scalar curvature.  In section 7 we show that the spectral triples we consider give rise to same Kasparov element. Then we also show that they have nontrivial Chern character.


\newsection{The Quantum Heisenberg Algebra}
Notation: for $x \in \IR$, $e(x)$ stands for $e^{2 \pi  ix}$
\bdfn
For any positive integer c let $S^c$ denote the space of  $ C^{\infty}$  functions
$\Phi : \IR \times \IT \times \IZ  \raro C$ such that \\
a) $\Phi (x+k,y,p)=e(ckpy) \Phi(x,y,p) $ for all $k \in \IZ$ \\
b) for every polynomial $P$ on $\IZ$ and every partial differential operator \\$ \widetilde{X}=\frac{\del^{m+n}}{\del x^m \del y^n} $ on $\IR \times \IT$ the function $ P(p)(\widetilde{X} \Phi) (x,y,p)$ is bounded on $ K \times \IZ$ for any compact subset $K$  of $\IR \times \IT$.\\
For each  $\hbar,\mu,\nu  \in \IR,\mu^2 + \nu^2 \ne 0$, let ${\cla}^{\infty}_{\hbar} $ denote $S^c$ with product and involution defined by
\bea
\label{1}
(\Phi \star \Psi)(x,y,p)= \sum_q \Phi(x-\hbar (q-p) \mu ,y-\hbar (q-p) \nu,q) \Psi (x- \hbar q \mu,y-\hbar q \nu,p-q) \eea
\bea \label{2} \Phi^*(x,y,p)= \bar{\Phi}(x,y,-p) \eea
$\pi : {\cla}^{\infty}_{\hbar}\raro \clb(L^2( \IR \times \IT \times \IZ)) $ given by
\bea \label{3}
(\pi (\Phi) \xi)(x,y,p)= \sum_q \Phi ( x - \hbar ( q-2p)\mu,y- \hbar (q-2p) \nu,q) \xi ( x,y,p-q)
\eea gives a faithful representation of the involutive  algebra $ \AAZ$. \\
${\cla}^{c,\hbar}_{\mu,\nu}=$ norm closure of $\pi( \AAZ )$ is called the Quantum Heisenberg Manifold.\\
$N_{\hbar}$= weak closure of $\pi ( \AAZ )$
\edfn
We will identify $\AAZ$ with $\pi( \AAZ) $ without any mention.\\
Since we are  going to work with fixed parameters $ c, \mu,\nu, \hbar $ we will drop them altogether and denote ${\cla}^{c,\hbar}_{\mu,\nu}$ simply by $\cla_\hbar$ here the subscript remains merely as a reminiscent of  Heisenberg only to distinguish it from a general algebra.\\


 {\bf Action of the heisenberg group:}
For $ \Phi \in S^c, (r,s,t) \in \IR^3 \equiv G$,  (as a topological space)
\bea
\label{4}
(L_{(r,s,t)} \phi )(x,y,p)=e(p(t+cs(x-r)))\phi(x-r,y-s,p) \eea
extends to an ergodic action of the Heisenberg group on $ {\cla}^{c,\hbar}_{\mu,\nu}$.\\


{\bf The Trace:}
 $\tau : \AAZ \raro \IC$, given by $\tau (\phi)= \int^1_0 \int_{\IT} \phi (x,y,0) dx dy $ extends to a faithful normal tracial state on $N_\hbar$.\\
$\tau$ is invariant under the Heisenberg group action. So, the group action can be lifted to $L^2( \AAZ)$.
We will denote the action at the Hilbert space level by the same symbol.
\bthm{{\bf (Weaver)}}
Let $ \clh = L^2(\IR \times \IT \times \IZ)$ and $V_f,W_k,X_r$ be the operators defined by \\
$(V_f \xi)(x,y,p)=f(x,y) \xi(x,y,p) $\\
$(W_k \xi) (x,y,p)= e(-ck ( p^2 \hbar \nu + py)) \xi(x+k,y,p)$\\
$(X_r \xi )( x,y,p)=\xi(x-2 \hbar r \mu,y-2 \hbar r \nu,p+r)$\\
Let $ T \in \clb (\clh)$. Then $T \in N_{\hbar}$  iff $T$ commutes with the operators $V_f,W_k,X_r$ for all $ f \in L^{\infty}       (\IR \times \IT ),k,r \in \IZ$.
\ethm


\blmma
Let $S^c_{\infty,\infty,1}=\{ \psi : \IR \times \IT \times \IZ \raro \IC |$ (i) $\psi$ is measurable,  \\  (ii) $\psi_n=\sup_{x \in \IR, y \in \IT} | \psi (x,y,p) | $ is an $ l_1$ sequence, \\ (iii) $ \psi (x+k,y,p)=e(ckyp)\psi(x,y,p) $ for all $ k \in \IZ \}$.\\
Then, for $\phi \in S^c_{\infty,\infty,1} $  $ \pi (\phi)$ defined by the same expression as in (\ref{3}) gives a bounded operator on $L^2(\IR \times \IT \times \IZ)$
\elmma
{\it Proof:} Let $\phi^{\prime} : \IZ \raro \IR_+ $ be $ \phi^{\prime} (n)= \sup_{x \in \IR,y \in \IT}|\phi(x,y,n)|$.\\ Then $| (\pi (\phi)\xi)(x,y,p)| \le ( \phi^{\prime} \star | \xi (x,y, . )|)(p)$,\\
where $\star$ denotes convolution on $\IZ$ and $|\xi(x,y, .)|$ is the function $ p \mapsto | \xi(x,y,p)|$.
By Young's inequality $\| (\pi(\phi) \xi) ( x,y, .) \|_{l_2}\le \| \phi^{\prime}  \star | \xi(x,y,.)| \|_{l_2} \le \|\phi^{\prime} \|_{l_1} \| \xi(x,y,.)\|_{l_2}$\\
Therefore , $ \| \pi ( \phi ) \| \le \| \phi \|_{\infty, \infty ,1 }$, where $\| \phi \|_{\infty,\infty,1}= \| \phi^{\prime} \|_{l_1}$ \qed.
\brmrk
i) product and involution defined by (\ref{1},  \ref {2}) turns $S^c_{\infty,\infty,1}$ into an involutive algebra.\\
ii) $ \phi \mapsto \| \phi \|_{\infty,\infty,1}$ is a $\star$-algebra norm.
\ermrk
\blmma
$\pi ( S^c_{\infty,\infty,1}  ) \subseteq N_\hbar.$
\elmma
{\it Proof:} Follows from Weaver's characterization of $N_\hbar$ \qed.


\bppsn
\label {26}
$L^2(\AAZ ,\tau)$ is unitarily equivalent with \\  $ L^2( \IT \times \IT \times \IZ)      \cong L^2( [0,1]\times [0,1] \times \IZ)$
\eppsn
{\it Proof:} For $ \phi \in S^c_{\infty,\infty,1}, $  $\Gamma \phi : \IR \times \IT \times \IZ \raro \IC$ is given by
\[
\Gamma \phi ( x,y,p)=\cases{e(-cxyp)\phi (x,y,p) &    for $ y <1 $\cr
\phi (x,y,p)&   for $ y=1$}.
\]
Then $\Gamma \phi ( x+k,y,p)=\Gamma \phi (x,y,p)$. So, $\Gamma \phi$ is a map from $\IT
\times \IT \times \IZ$ to $\IC$.\\
$$ \tau (\phi^{*} \star \phi ) = \int^1_0 \int_{\IT} \sum_q {| \phi(x- \hbar q \mu,y- \hbar q \nu,-q)|}^2dxdy
 =\int^1_0 \int_{\IT} \sum_q { | \phi (x,y,q) |}^2dxdy $$
Since, $| \phi ( x+k,y,p)|=|\phi (x,y,p) |$ for all $x \in \IR,y \in \IT, k,p \in \IZ$.\\
Therefore  $ \tau (\phi^{*} \star \phi)= {\| \Gamma \phi \|}^2$, i.e, $\Gamma : L^2(\AAZ,\tau) \raro
L^2( \IT \times \IT \times \IZ) $ is an isometry. To see $ \Gamma $ is an unitary observe,\\
(i) $N_{\hbar} \subseteq L^2(\AAZ,\tau)$, since $\tau$ is normal.\\

(ii)
\[
 \phi_{m,n,k}= \cases { e(cxyp) e(mx+ny) \delta_{kp},  &  for   $0 \le y \le 1$  \cr
   \delta_{kp} e(mx)  & for   $y=1$ }
\]
 is an element of $S^c_{\infty, \infty,1} \subseteq N_{\hbar}$\\
(iii) ${\{ \Gamma \phi_{m,n,k} \}}_{m,n,k\in \IZ}$ is an orthonormal basis in $ L^2({\IT}^2 \times \IZ).$ \qed
\brmrk
$\phi \mapsto \phi |_{[0,1]\times \IT \times \IZ} $ gives an unitary isomorphism.
\ermrk
\bcrlre
\label {28}
Let $M_{yp} $ be the multiplication operator on $\clh= L^2(\IT \times \IT \times \IZ)$. If we consider $\AAZ$ as a subalgebra of  $\clb (\clh) $ by the left regular representation then,
$[M_{yp}, \AAZ ] \subseteq  \clb (\clh)$
\ecrlre
{\it Proof:} Note for $ \phi \in \AAZ$, $(M_{yp} \phi)(x,y,p)=yp \phi(x,y,p)$ gives an element in $S^c_{\infty,\infty,1}$ hence a bounded operator.
\bean
\lefteqn{[M_{yp},\phi] \psi (x,y,p)}\\
&=& \sum_q (yp-(y-\hbar q \nu)(p-q))\phi(x-\hbar (q-p)\mu, y- \hbar (q-p) \nu,q) \\
& & \times \psi ( x- \hbar q \mu , y - \hbar q \nu,p-q)\\
&=& \sum_qq (y-\hbar (q-p) \nu) \phi ( x- \hbar (q-p) \mu, y- \hbar (q-p) \nu,q) \\
& & \times\psi ( x-\hbar q \mu, y- \hbar q \nu,p-q )\\
&=& ( M_{yp} (\phi) \star \psi )(x,y,p)\\
\eean
for $ \psi \in \AAZ$.
This completes the proof. \qed


\newsection {A class of spectral triples}


Let $ (\cla ,G ,\al)$ be a $ C^*$ dynamical system with $G$ an n dimensional Lie group, and $\tau$ a $G$-invariant trace on $\cla$. Let $\cla^{\infty}$ be the space of smooth vectors, \\$ {\cal K}= L^2(\cla,\tau) \otimes {\IC}^N$ where $N=2^{\lfloor n/2 \rfloor}$. Fix any basis $X_1,X_2,\ldots X_n$ of $L(G)$ the Lie algebra of G. Since G acts as a strongly continuous unitary group on $ \clh = L^2(\cla,\tau)$ we can form selfadjoint operators $d_{X_i}$ on $\clh$. $D: {\cal K}\raro {\cal K}$ is given by
$D= \sum_i d_{X_i}\otimes \gamma_i$, where $ \gamma_1,\ldots \gamma_n$ are selfadjoint matrices in $M_N(\IC)$ such that $ \gamma_i \gamma_j + \gamma_j \gamma_i= 2 \dlt_{ij}$  along with $\cla^{\infty}$ and ${\cal K}$ should produce a spectral triple. For such a $D$, $[D,\cla^{\infty}] \subseteq \cla^{\infty}\otimes M_N(\IC)$.


\bppsn
\label {31}
For the quantum Heisenberg manifold, if we identify the Lie algebra of Heisenberg group with the Lie algebra of upper triangular matrices, then D as described becomes a selfadjoint operator with compact resolvent for the following choice: $$X_1=\left(\matrix { 0&1&0 \cr 0&0&0 \cr 0& 0&0\cr}\right),
X_2=\left(\matrix { 0&0&0 \cr 0&0&1\cr 0& 0&0\cr}\right),X_3=\left(\matrix { 0&0&c \al \cr 0&0&0 \cr 0& 0&0\cr}\right)$$
where $ \al \in \IR$ is greater than one.
\eppsn
{\it Proof:}
Let $\cld (D)=\{ f \in L^2([0,1]\times [0,1]\times \IZ) | f(x,0,p)=f(x,1,p) , $\\ $f(1,y,p)=e(cpy)f(0,y,p),pf,\frac{\del f}{\del x},\frac{\del f}{\del y} \in L^2 \} \otimes \IC^2$\\
$D(f \otimes u)= \sum^3_{j=1} i d_j(f) \otimes \sigma_j(u)$, where
\bean
\label{8}
id_1(f)&=& -i \frac {\del f } {\del x} \cr
i d_2 (f) &=& -2 \pi c p x f(x,y,p) - i \frac { \del f } {\del y} \cr
id_3 (f) &=& -2 \pi p c \al f(x,y,p) \cr
\eean
and $\sigma_j$'s are the spin matrices.\\
Let $\eta :L^2([0,1]\times [0,1]\times \IZ) \raro L^2([0,1]\times [0,1] \times \IZ)$ be the unitary given by
\[
\eta (f)(x,y,p)=\cases{e(-cxyp)f(x,y,p)& for $y <1$ \cr
f(x,y,p)& for $y=1$}
\]
Let $\cld (D^{\prime})= (\eta \otimes I_2 ) \cld (D)$, and $D^{\prime}=(\eta \times I_2) D {(\eta \otimes I_2)}^{-1}$.\\Then $\cld ( D^{\prime}) =\{ f\in L^2([0,1] \times [0,1] \times \IZ) | f(x,0,p)=f(x,1,p),$\\ $f(0,y,p)=f(1,y,p),
\frac{\del f}{\del x}, \frac {\del f}{\del y},pf \in L^2 \} \otimes \IC^2$ and\\
$D^{\prime}( f \otimes u)= \sum^{3}_{j=1}id^{\prime}_j (f) \otimes \sigma_j (u)$ where,
\bean
d^{\prime}_1(f)(x,y,p)&=& -2 \pi i c y p f(x,y,p) - \frac {\del f}{\del x} (x,y,p) \cr
d^{\prime}_2(f)(x,y,p) &=& -\frac{\del f}{\del y}(x,y,p)  \cr
d^{\prime}_3 (f)(x,y,p)& =& 2 \pi i p c \al f(x,y,p)  \cr
\eean
Note, on $\cld (D^{\prime}),D^{\prime}= T+S$ where, $ \cld  (T)= \cld (D^{\prime})\subseteq \cld (S)$\\
$T=-i\frac{\del}{\del x}\otimes \sigma_1 -i \frac{\del}{\del y} \otimes \sigma_2 -2 \pi c \al M_p \otimes \sigma_3$\\
$S=2 \pi c M_{yp}  \otimes \sigma_1$ are selfadjoint operators on their respective domains.  Also note $T$ has compact resolvents. Our conclusion  follows from the Rellich lemma since $S$ is relatively bounded with respect to $T$ with relative bound less than $\frac {1}{\al} <1$. \qed


\bthm
Let $ \clh= L^2(\AAZ,\tau)\otimes \IC^2,$. $\AAZ$  with its diagonal action becomes a subalgebra of $\clb (\clh)$.
$(\AAZ, \clh,D) $ is an odd spectral triple of dimension 3.
\ethm
{\it Proof:} $(\AAZ,\clh,D)$ is a spectral triple follows from the previous proposition and the remark preceeding that. We only have to show ${|D|}^{-3} \in {\cll}^{(1,\infty)}$, the ideal of Dixmier traceable operators. For that observe:\\
(i) Since $T$ is the dirac operator on ${\IT}^3$, $\mu_n(T^{-1} |_{{ker T}^\perp})=O(1/n^{1/3})$, $\mu_n$ stands for the nth singular value.\\
(ii)  S is relatively bounded with relative bound less than $\frac {1}{\al} < 1$, hence
$\|S{(T+i)}^{-1}\| \le \frac{1}{\al} $ and  $ \| {(1+S{(T+i)}^{-1})}^{-1}\| \le \frac{\al}{\al-1}$ \\
(iii) $\mu_n(AB) \le \mu_n(A) \|B\|$, for bounded operators A,B.\\
Applying (i),(ii),(iii) to ${(D^{\prime} +i)}^{-1}= {(T+i)}^{-1} { (1+S{(T+i)}^{-1})}^{-1} $ we get
the desired conclusion for $D^{\prime}$ and hence for D. \qed
\bcrlre
\label {33}
Let $T,S,D,D^\prime$ be as in the proof of proposition (\ref {31}).\\$ A=  {(\eta \otimes I_2)}^{-1} T
(\eta \otimes I_2)$  Then $ (\AAZ, \clh,A) $ is an odd spectral triple of dimension 3
\ecrlre
{\it Proof:}We only have to show $[A, \AAZ] \subseteq \clb ( \clh) $.\\
Let $ B= {(\eta \otimes I_2)}^{-1} S
(\eta \otimes I_2)$. Then since $ \eta \otimes I_2$ commutes with $S$, B=S. By corollary (\ref {28}),
$[B,\AAZ] \subseteq  \clb ( \clh) $. Now the previous theorem along with $D=A+B$ completes the proof. \qed\\
\brmrk
One can similarly show $ (\AAZ, \clh,A_t) $ forms an odd spectral triple of dimension 3, for $t \in [0,1]$. Here $ A_t$ stands for $A_t= A+ tB$
\ermrk
\brmrk
D,A constructed above depends on $\al$.
\ermrk


\bppsn
\label {36}
The positive linear functional on $\cla_{\hbar} \otimes M_2(\IC)$ given by $ \int : a \mapsto {tr}_\omega a {|D|}^{-3}$ is nothing but $\frac{1}{2}({tr}_\omega {|D|}^{-3}) \tau \otimes tr$.
\eppsn
{\it Proof :}
\bean
D^2  & = & -\left( \matrix {d_1^2 + d_2^2 + {(d_3+\frac{1}{2 \al})}^2-\frac{1} {4{\al^2}} &0 \cr
0 &  d_1^2 + d_2^2 + {(d_3-\frac{1}{2 \al})}^2-\frac{1} {4{\al^2}} \cr }\right)   \cr
     & =  & \left( \matrix {X_1&0 \cr 0 & X_2 \cr} \right),  {\mbox say.} \cr
\eean
  It is easily seen that  \\ (i)  compactness of  resolvents of $D^2$ implies that for $X_1,X_2$ \\ (ii) eigenvalues of $X_1,X_2$ have similar asymptotic behaviour.   \\
Therefore $ X_1^{-3/2}, X_2^{-3/2} \in {\cll}^{(1,\infty)}$ and $ tr_{\omega}  a  X_1^{-3/2}=tr_{\omega} a  X_2^{-3/2}$ for any $ a \in \clb (L^2(\cla_\hbar))$\\
Consider the unitary group on $ \clh \cong L^2([0,1]\times \IT \times \IZ) \otimes \IC^2$ given by \\
$U_t ( x \otimes y \otimes e_p \otimes z)=e(pt) ( x \otimes y \otimes e_p \otimes z).$\\
Then $ U_t D= D U_t$ and
$$\int A = tr_\omega U_t A U_t^* {|D|}^{-3}=  tr_\omega
\int_0^1 U_t A U_t^* {|D|}^{-3} dt = \int {(A)}_0 $$
Here  $$ A= \left(  \matrix { \psi_{11} & \psi_{12} \cr \psi_{21} & \psi_{22}  \cr } \right) \mapsto {(A)}_0=
\left( \matrix { (\psi_{11})_0 & (\psi_{12})_0  \cr      (\psi_{21})_0  &     (\psi_{22})_0 \cr } \right) $$ is the CPmap explicitly given for $ \psi \in S^c$ by $ (\psi)_0(x,y,p)= {\delta}_{p0} \psi(x,y,p)$\\
Since $\left( \matrix { 1 & 0 \cr 0 & 0 \cr } \right) $ commutes with ${|D|}^{-3}$, we get
\bean
\int A&= &tr_\omega {(a_{11})}_0 {X_1}^{-3/2} +tr_\omega {(a_{22})}_0 {X_2}^{-3/2}\cr
&=& tr_\omega({(a_{11})}_0 +{(a_{22})}_0) {X_1}^{-3/2} \cr
\eean
Consider the homomorphism $ \Phi : C(\IT^2) \raro \cla_\hbar$ given by \\ $ \Phi(f)(x,y,p)=\delta_{p0} f(x,y)$.
Now by riesz representation theorem for \\ $ \int \circ (\Phi \otimes I_2) : C(\IT^2) \raro \IC$, we get a measure $\lmd $ on $\IT^2$ such that \\ $ tr_\omega 2 {(\psi)}_0 {X_1}^{-3/2}=\int {(\psi)}_0(x,y,0) d \lmd$
implying
\bea
\label{6}
\int A= \frac{1}{2} \int ({(a_{11})}_0 +{(a_{22})}_0)d \lmd
\eea
In the next lemma we show $ \lmd $ is proportional to lebesgue measure. That will prove $\int \propto \tau \otimes tr$, and the proportionality constant is obtained by evaluating both sides on $I$. \qed
\blmma
If $\{ 1,\hbar \mu , \hbar \nu \}$ is rationally independent then $\lmd$ as obtained in the previous proposition is proportional to lebesgue measure.
\elmma
{\it Proof:} It is known \cite {VAR} \cite {CGS} that for a spectral triple $(\cla, \clh, \cld)$ with ${|\cld|}^{-p} \in {\cll}^{(1,\infty)} $ for some p, $ a \mapsto tr_\omega a {|\cld|}^{-p}$ is a trace on the algebra. This along with (\ref{6}) gives
\[
\label{7}
\int (\phi \star \psi)(x,y,0) d\lmd (x,y)= \int (\psi \star \phi)(x,y,0) d \lmd (x,y) ,
\forall \phi,\psi \in S^c
\]
Taking $\phi (x,y,p)=e(c [x]yp)f(x-[x]) g(y) \delta_{1p}$ where $g: \IT \raro \IC,f:[0,1] \raro \IC$ are smooth functions with $supp(f) \subseteq [\epsilon,1-\epsilon]$ for some $\epsilon >0$  and $ \psi = {\phi}^* $ we get
from (\ref{7})
\[
\label {8}
\int{ | \phi (x+\hbar \mu,y+\hbar \nu ,1)|}^2 d \lmd (x,y)
=\int{|\phi \circ \gamma ( x+\hbar \mu, y+ \hbar \nu, 1)|}^2 \lmd(x,y)
\]
where $ \gamma :\IT^2 \raro \IT^2 $ is given by $ \gamma (x,y)=(x-2 \hbar \mu,y-2 \hbar \nu)$.
The hypothesis of linear independence of $(1,\hbar \mu,\hbar \nu)$ over the rationals implies that $\gamma$-orbits are dense. This along with (\ref {8}) proves the lemma. \qed
\brmrk
In the rest of the paper $\int$ will denote $ \frac{1}{2} \tau \otimes tr$.
\ermrk


\newsection{Space of forms }
\blmma
Let $\cla$ be a dense subalgebra of a unital C*algebra $\bar{\cla}$ closed under holomorphic function calculus, then $\cla$ is simple provided $\bar {\cla} $ is so.
\elmma
{\it Proof:}
Let $ J \subseteq \cla$ be an ideal. Then $\bar{J}=\bar{\cla}$, since $\bar{\cla}$ is simple. There exists $x \in J $ such that $\| x-I\| <1$. Then $x^{-1} \in \bar{\cla} $, hence in $\cla$ because $\cla$ is closed under holomorphic function calculus. Therefore $1=x x^{-1} \in J$. \qed
\brmrk
$\AAZ$ is simple because $\cla_\hbar $ is so.
\ermrk


\bdfn {\bf (Connes) }
Let $(\cla ,\clh ,D)$ be a spectral triple.\\
$\Omega^k(\cla)=\{ \sum_i=1^N a_0^{i} \delta a_1^i \ldots \delta a_k^i | n \in \IN, a_j^i \in \cla \} { } \Omega^\bullet(\cla)=\oplus_0^{\infty} \Omega^k(\cla)$
is the unital graded algebra of universal forms. Here $\delta $ is an abstract linear operator with $\delta^2=0,\delta(ab)=\delta(a)b+a\delta(b)$. $\Omega^\bullet (\cla)$ becomes a *algebra under the involution ${(\delta a)}^*=-\delta (a^* ) \forall a \in \cla$. Let $ \pi: \Omega^\bullet ( \cla)  \raro \clb (\clh)$ be the $\star$-representation given by $ \pi(a)=a,\pi (\delta a)=[D,a]$
Let $J_k= ker \pi |_{\Omega^k(\cla)}$ The unital graded differential $\star$-algebra of differential forms
$\Omega^\bullet_D(\cla)$ is defined by
\[
\Omega_D^\bullet ( \cla ) = \oplus_0^\infty \Omega_D^k(\cla), \Omega_D^k(\cla)= \Omega^k(\cla)/(J_k+\delta J_{k-1} \cong \pi (\Omega^k(\cla))/\pi ( \delta J_{k-1})
\]
\edfn
Notation:--(i)Let $\phi \in S^c$, then $[D,\phi]=\sum \delta_i (\phi) \otimes \sigma_i$ where
$\delta_j(\phi)= i d_j(\phi)$ (see proof of proposition \ref {31} for $d_j$ ) but looked upon as derivation on $\AAZ$.\\
Note: $[\delta_1,\delta_3]=[\delta_2,\delta_3]=0, [\delta_1,\delta_2]=\delta_3$\\
(ii) $ \phi_{m,n}(x,y,p)=e(mx+ny)\delta_{p0}$
\blmma \label{tl1}
Let $\cla$ be a unital simple algebra, $M \subseteq \underbrace {\cla \oplus \ldots \oplus\cla}_{n times}$ a sub$\cla$-$\cla$ bimodule. Suppose $\exists a_{ij},1 \le n, 1 \le j \le i$ such that \\
(i) $a_{ii} \ne 0$, (ii) $b_i=(a_{i1},\ldots ,a_{ii},0,\ldots,0) \in M$\\
Then $M \cong \underbrace {\cla \oplus \ldots \oplus\cla}_{n times}$ as an $\cla$-$\cla$ bimodule.
\elmma
{\it Proof:} By induction on n,\\
For n=1, $0 \ne M$ is an ideal in $\cla$, hence $M=\cla$.\\
Let $ \pi : M \raro \cla$ be $ \pi (a_1,\ldots , a_n)=a_n$.\\
Then by hypothesis $\pi (M) $ is a nontrivial ideal in $\cla$ hence equals $\cla$. So, we have a split short exact sequence $$
0 \raro \ker(\pi) \raro M \raro \cla \raro 0$$
Therefore $M=ker(\pi) \oplus Im \pi = ker(\pi) \oplus \cla =\underbrace {\cla \oplus \ldots \oplus\cla}_{n times}$. In the last equality we have used induction hypothesis for $ker (\pi)$. \qed


\bppsn
\bean
 (i) \Omega^1_D ( \AAZ) & =  & \{ \sum a_i \otimes \sigma_i| a_i \in \AAZ , \sigma_i 's  \mbox { are spin  matrices  } \}   \cr
    & = & \AAZ \oplus \AAZ \oplus \AAZ  \cr
\eean
(ii)
$ \pi (\Omega^k (\AAZ))  = \AAZ \otimes M_2(\IC)= \AAZ \oplus \AAZ \oplus \AAZ \oplus \AAZ$
\eppsn
{\it Proof:}$\Omega^1_D(\AAZ) = \pi ( \Omega^1(\AAZ)) \subseteq $ R.H.S.\\
Let $\phi_{m,n}(x,y,p)=\delta_{p0} e(mx+ny)$ and $\phi \in S^c$ be such that $ \phi (x,y,p)=\delta_{p1} \phi(x,y,p)$. Then applying the previous lemma to $[D,\phi_{01}],[D,\phi_{10}],[D,\phi] \in \pi (\Omega^1 (\cla))$ we get the result. \\
(ii)(i) along with $ \Omega^k (\AAZ)=  {\underbrace  {\Omega^1 (\AAZ) \otimes_{\AAZ}  \ldots \otimes_{\AAZ} \Omega^1 (\AAZ)}_{k times }}$ proves the result. \qed


\bppsn
(i) $\pi ( \delta  J_1)=\AAZ $\\
(ii) $ \Omega^2_D ( \AAZ) = \AAZ \oplus \AAZ \oplus \AAZ$
\eppsn
{\it Proof:}(i) Let $\omega =\sum a_i \delta (b_i) \in J_1$. Then $ \pi (\omega)= \sum a_i \delta_j (b_i) \sigma_j =0$ gives $ \sum a_i \delta_j(b_i)=0, \forall j$.
\bea
\label{9}
\pi(\delta \omega)&=& \sum_i (\sum_j \delta_j ( a_i) \sigma_j )( \sum_k \delta_k (b_i ) \sigma_k)   \cr
&=& \sum_i ( \sum_j \delta_j (a_i) \delta_j (b_i) )\otimes I_2 +\sum_i (   \sum_{j<k}
(\delta_j(a_i) \delta_k(b_i)- \delta_k(a_i) \delta_j(b_i)) \sigma_j \sigma_k)\cr
\eea
\bea
\label {10}
\sum_i [\delta_j,\delta_k](a_ib_i)&=& \sum_i  \delta_j ( \delta_k (a_i) b_i) - \delta_k(\delta_j(a_i) b_i) [ {Since}  \sum a_i \delta_j (b_i)=0, \forall j \cr
&=& \sum_i [ \delta_j , \delta_k ] (a_i) b_i +\sum_i (\delta_k(a_i) \delta_j(b_i)-\delta_j (a_i) \delta_k(b_i)) \cr
\eea
Also note
\bea
\label{11}
\sum_i [\delta_j,\delta_k](a_ib_i)&=&\sum_i [\delta_j,\delta_k](a_i)b_i+\sum_i a_i[\delta_j,\delta_k](b_i) \cr
&=& \sum_i \sum_i [\delta_j,\delta_k](a_i)b_i \cr
\eea
Comparing rhs of (\ref{10},\ref {11}) we see second term on the rhs of (\ref{9}) vanishes proving $\pi (\delta J_1) \subseteq \AAZ$. \\
To see actually equality holds note,\\
$\omega = 2 \phi_{0,2} \delta (\phi_{0,1}) - \phi_{0,1} \delta (\phi_{0,2}) \in J_1 \pi(\delta \omega)= 2 \phi_{0,3} \otimes I_2 \ne 0$.  An application of lemma~\ref{tl1} proves (i).\\
(ii)Let $ \phi \in S^c $ be such that $ \phi(x,y,p)=\delta_{1p} \phi(x,y,p)$ \\
$\omega_1= \delta (\phi_{1,0}) \delta (\phi_{0,1})
\omega_2= \delta (\phi_{1,0}) \delta (\phi)
\omega_1= \delta (\phi_{0,1}) \delta (\phi) $
Now lemma~\ref{tl1} together with (i) implies the result. \qed


\blmma
$\pi (\delta J_2)= \{\sum a_j \otimes \sigma_j| a_j \in \AAZ \}= \AAZ \oplus \AAZ \oplus \AAZ$
\elmma
{\it Proof:} Let $\omega = \sum a_i \delta(b_i) \delta(c_i)  \in J_2$
\bean
0=\pi(\omega)&=& \sum a_i(\sum_j \delta_j (b_i) \sigma_j)(\sum_j \delta_k (c_i) \sigma_k) \\
&=&  a_i \delta_j(b_i) \delta_j(c_i) + \sum_{j<k} a_i (\delta_j(b_i) \delta_k (c_i)- \delta_k(b_i) \delta_j(c_i)) \sigma_j \sigma_k \\
\eean
Comparing the coefficients of the various spin matrices we get
\bea
\label{12}
\sum a_i \delta_j(b_i) \delta_j (c_i) &=&0 \\
\sum a_i (\delta_j(b_i) \delta_k(c_i)-\delta_k(b_i)\delta_j(c_i))&=&0, \forall j \ne k \\
\eea
from (\ref{12}),
\bean
0 &=& \sum \delta_1 ( a_i(\delta_2(b_i)\delta_3(c_i)-\delta_3(b_i) \delta_2( c_i)) \cr
&=& \sum \delta_1(a_i) ( \delta_2(b_i)\delta_3(c_i)-\delta_3(b_i) \delta_2( c_i)) \cr
& &+ \sum a_i \delta_1 (\delta_2(b_i)\delta_3(c_i)-\delta_3(b_i) \delta_2( c_i)) \cr
\eean
Therefore,
\bea
\sum \delta_1(a_i)(\delta_2(b_i)\delta_3(c_i)-\delta_3(b_i) \delta_2( c_i)) &=&
-\sum a_i \delta_1(\delta_2(b_i)\delta_3(c_i)-\delta_3(b_i) \delta_2( c_i)) \cr
\eea
Similarly we get two more equalities.
Let $A=$ be coefficient of $I_2$ in $\pi (\delta \omega)$. Then
\bean
\sqrt{(-1)}A&=& \sum \delta_1(a_i)(\delta_2(b_i)\delta_3(c_i)-\delta_3(b_i) \delta_2( c_i)) \cr
& &+\sum \delta_2(a_i)(\delta_3(b_i)\delta_1(c_i)-\delta_1(b_i) \delta_3( c_i))   \cr
& &+ \sum \delta_3(a_i)(\delta_1(b_i)\delta_2(c_i)-\delta_2(b_i) \delta_1( c_i)) \cr
&=& -(\sum  a_i \delta_1 (\delta_2 (b_i) \delta_3 (c_i) - \delta_3(b_i) \delta_2(c_i)) \cr
& &+\sum  a_i \delta_2 (\delta_3(b_i) \delta_1 (c_i) - \delta_1(b_i) \delta_3(c_i)) \cr
& & +\sum  a_i \delta_1 (\delta_1 (b_i) \delta_2 (c_i) - \delta_2(b_i) \delta_1(c_i))) \cr
&=& -(\sum a_i ( [\delta_1,\delta_2](b_i) \delta_3(c_i) +\delta_2 (b_i) [\delta_1,\delta_3](c_i)) \cr
& &+\sum a_i([\delta_3,\delta_1](b_i) \delta_2(c_i) +\delta_3 (b_i) [\delta_1,\delta_3](c_i)) \cr
& &+\sum a_i( [\delta_2,\delta_3](b_i) \delta_1(c_i) +\delta_1 (b_i) [\delta_3,\delta_2](c_i))) \cr
      & =& 0 \cr
\eean
Here second equality follows from (4.7) and the last equality follows from (4.5) since $\delta_j$'s form a lie algebra.
This shows
\[
\label{13}
 \pi (\delta J_2) \subseteq  \{ \sum_{j=1}^{3} a_j \sigma_j | a_j \in \AAZ \}\cong \AAZ \oplus \AAZ \oplus  \AAZ
\]
Let $\phi \in S^c$ be such that $ \phi(x,y,p)=\delta_{1p} \phi(x,y,p)$.
Then,
\bean
\omega_1 &=& 2 \phi_{0,2} \delta (\phi_{0,1}) \delta ( \phi_{0,1})-  \phi_{0,1} \delta (\phi_{0,2}) \delta ( \phi_{0,1}) \in J_2 \cr
\omega_2 &=& 2\phi_{2,0} \delta (\phi_{1,0}) \delta ( \phi_{1,0})-\phi_{1,0} \delta (\phi_{2,0}) \delta ( \phi_{1,0}) \in J_2 \cr
\omega_3 &=& \phi_{0,2} \delta (\phi_{0,1}) \delta ( \phi)-\phi_{0,1} \delta (\phi_{0,2}) \delta ( \phi)
 \in J_2 \cr
\eean
satisfies,
\bean \pi (\delta \omega_1)  &=&2 \phi_{0,4} \sigma_2 \cr
\pi (\delta \omega_2 ) &=&  2 \phi_{4,0} \sigma_1 \cr
\pi (\delta \omega_3) &=&  2 \phi_{0,3} \delta_1(\phi) \sigma_1+ 2 \phi_{0,3} \delta_2(\phi) \sigma_2
+ 2 \phi_{0,3} \delta_3(\phi) \sigma_3 \cr
\eean
Therefore by Lemma~\ref {tl1} we get equality in (\ref {13}).
\qed.
\bcrlre
$\Omega^3_D (\AAZ) =\AAZ$
\ecrlre
{\it Proof:} Immediate from the previous lemma and proposition 3.5(ii). \qed


\blmma
(i)$ \Omega^4_D( \AAZ )= 0 $\\
(ii) $ \Omega^k_D( \AAZ )= 0, \forall k > 4 $
\elmma
{\it Proof:} (i) It suffices to show $ \pi ( \delta J_3)= \AAZ \oplus  \AAZ \oplus  \AAZ \oplus
\AAZ $. \\
For that note,
\bean
\omega_1 &=& 2 \phi_{0,2} \delta(\phi_{0,1})     \delta(\phi_{0,1})   \delta(\phi_{0,1})
-\phi_{0,1}  \delta(\phi_{0,2})   \delta(\phi_{0,1})  \delta(\phi_{0,1})  \in J_3 \cr
\omega_2 &=& 2 \phi_{0,2} \delta(\phi_{0,1})     \delta(\phi_{0,1})   \delta(\phi_{0,1})
-\phi_{0,1}  \delta(\phi_{0,2})   \delta(\phi_{0,1})  \delta(\phi_{0,1})  \in J_3 \cr
\omega_3 &=& 2 \phi_{0,2} \delta(\phi_{0,1})     \delta(\phi_{0,1})   \delta(\phi)
-\phi_{0,1}  \delta(\phi_{0,2})   \delta(\phi_{0,1})  \delta(\phi)  \in J_3 \cr
\omega_4 &=& 2 \phi_{0,2} \delta(\phi_{0,1})     \delta(\phi_{1,0})   \delta(\phi)
-\phi_{0,1}  \delta(\phi_{0,2})   \delta(\phi_{1,0})  \delta(\phi)  \in J_3 \cr
\eean
satisfies
\bean
  \pi (\delta \omega_1)&=& 2 \phi_{0,5} \otimes I_2 \cr
 \pi (\delta \omega_2)&=& 2 \phi_{1,4} \sigma_2 \sigma_1 \cr
\pi (\delta \omega_3)&=&2 \phi_{0,4} \delta_2 (\phi) \otimes I_2  +2 \phi_{0,4} \delta_1 (\phi) \sigma_2 \sigma_1  +2 \phi_{0,4} \delta_3 (\phi) \sigma_2 \sigma_3 \cr
\pi (\delta \omega_4)&=& 2 \phi_{1,3} \delta_1 (\phi) I_2   +2 \phi_{1,3} \delta_2 (\phi) \sigma_1 \sigma_2
+2 \phi_{1,3} \delta_3 (\phi) \sigma_1  \sigma_3  \cr
\eean
Now an application of Lemma~\ref{tl1} completes the proof. \\
(ii) The same argument as in (i) does the job with the following choice,\\
$ \omega_i^{\prime} = \omega_i \underbrace {\delta (\phi_{0,1}) \ldots \delta(\phi_{0,1})}_{(k-4) times},i=1,\cdot  4$ \qed


\newsection { Connections:-- torsinless/unitary}

\bdfn
{\bf {\cite {FGR}}}
(i)$ \int$ determines a semi-definite sesquilinear form on $\Omega^\bullet(\AAZ)$ by setting
$$ (\omega,\eta)=  \int \pi (\omega ) {\pi (\eta)}^* \forall \omega,\eta \in \Omega^\bullet( \AAZ) $$
(ii)Let $$ K_k=\{ \omega \in \Omega^k (\AAZ)_ ( \omega,\omega)=0 \}, K=\oplus^\infty_{k=0} K_k $$
$K, K+ \delta K $ are two sided *-ideals, the later is closed under differential.
$${\widetilde {\Omega}}^\bullet ( \AAZ) =\oplus^\infty_{k=0}  {\widetilde {\Omega}}^k ( \AAZ),
{\widetilde {\Omega}}^k ( \AAZ) = \Omega^k(\AAZ) /K_k$$
(iii)${\widetilde {\clh}}^k $ denotes the Hilbert space completion of  ${\widetilde {\Omega}}^k ( \AAZ) $
with respect to the scalar product. ${\widetilde {\clh}}^\bullet= \oplus^\infty_{k=0} {\widetilde {\clh}}^k $,
${\widetilde {\clh}}^k$ is to be interpreted as the space of square-integrable k-forms.\\
(iv)The algebra multiplication of  ${\Omega}^\bullet ( \AAZ) $ descends to a linear map \\$ m:
{\widetilde {\Omega}}^\bullet ( \AAZ) \otimes_{\AAZ}  {\widetilde {\Omega}}^\bullet ( \AAZ)
\raro   {\widetilde {\Omega}}^\bullet ( \AAZ) $.\\
(v) The   unital graded differential *-algebra of square-integrable differential  forms is defined by
$$ {\widetilde {\Omega}}^\bullet_D( \AAZ) =\oplus^\infty_{k=0} {\widetilde {\Omega}}^k_D ( \AAZ)
,   {\widetilde {\Omega}}^k_D ( \AAZ)=
{\widetilde {\Omega}}^k ( \AAZ)   /K_k + \delta K_{k-1}$$
(vi) $ \delta :{\Omega}^{\bullet +1} ( \AAZ)  \raro
{\Omega}^{\bullet +1}  ( \AAZ) $ descends to a linear map \\ $ \delta :
{\widetilde {\Omega}}^{\bullet}_D ( \AAZ) \raro  {\widetilde {\Omega}}^{\bullet  +1}_D ( \AAZ)$\\
(vii)  A {\bf connection} $\nabla$on a finitely generated projective $\AAZ$ module ${\cal E}$is a $ \IC$ linear map
$$ \nabla : {\widetilde {\Omega}}^{\bullet}_D ( \AAZ) \otimes {\cal E} \raro  {\widetilde {\Omega}}^{\bullet+1}_D ( \AAZ)  \otimes {\cal E}  $$
such that $\nabla (\omega s)= \delta(\omega)s+{(-1)}^k \omega \nabla (s) $ for all $ \omega \in
{\widetilde {\Omega}}^{\bullet}_D ( \AAZ)$ and all $ s \in {\widetilde {\Omega}}^{\bullet}_D ( \AAZ) \otimes {\cal E} $\\
(viii) The {\bf curvature} of a connection $\nabla$ on ${\cal E} $ is given by
$$ R(\nabla) =- {\nabla}^2 : {\cal E} \raro  {\widetilde {\Omega}}^k_D ( \AAZ) \otimes_{\AAZ} {\cal E}$$
\edfn
\brmrk
\label {52}
$ \omega \in  {\widetilde {\Omega}}^k ( \AAZ)$ determines two  operators \\$ m_L(\omega)
,  m_R(\omega)
: {\widetilde {\Omega}}^n ( \AAZ) \raro {\widetilde {\Omega}}^{n+k} ( \AAZ) $ given by $
m_L(\omega) (\eta)=m( \omega \otimes \eta)$,  \\  $m_R(\omega) (\eta)=m( \eta \otimes \omega)$. These operators extend to bounded linear operators $ m_L(\omega)
,m_R(\omega) : {\widetilde {\clh}}^n  \raro {\widetilde {\clh}}^{n+k}$ for all n.
\ermrk


\bppsn
(i)${\widetilde {\Omega}}^k ( \AAZ )= \AAZ \otimes M_2(\IC) \cong \AAZ \oplus  \AAZ \oplus
\AAZ \oplus  \AAZ$ \\
(iii) $ {\widetilde {\clh}}^k= L^2( \AAZ, \tau) \otimes \IC^4 $\\
(ii) ${\widetilde {\Omega}}^{k}_D ( \AAZ) =  {\Omega}^{k}_D ( \AAZ) $
\eppsn
{\it Proof:}(i) By the faithfulness of the linear functional $ A \mapsto \int A $ defined on
$\pi (\Omega^{\bullet}(\AAZ) )= \AAZ \otimes M_2(\IC)$ we get $J_k=K_k$. \\
hence ${\widetilde {\Omega}}^k ( \AAZ )= {\Omega}^{k} ( \AAZ)/ker(\pi) \cong \pi({\Omega}^{k} ( \AAZ))=\AAZ \otimes M_2(\IC)$\\
(ii) Follows from (i) and proposition \ref {36}\\
(iii) In (i) we have already seen $J_k=K_k$. That gives the result.  \qed
\brmrk
Since ${\widetilde {\Omega}}^{1}_D ( \AAZ)$ is free with 3 generators, we can and will identify
${\widetilde {\Omega}}^{1}_D ( \AAZ) \otimes_{\AAZ} {\widetilde {\Omega}}^{1}_D ( \AAZ)$
with $\AAZ \otimes M_3(\IC)$ and a connection $\nabla$ is specified by its value on the generators.
\ermrk


\bdfn
A connection $\nabla :{\widetilde {\Omega}}^{1}_D ( \AAZ)
\raro {\widetilde {\Omega}}^{1}_D ( \AAZ)    \otimes_{\AAZ}
 {\widetilde {\Omega}}^{1}_D ( \AAZ)$ is called torsionless if $ T(\nabla)= \delta -m\circ \nabla : {\widetilde {\Omega}}^{1}_D ( \AAZ) \raro {\widetilde {\Omega}}^{2}_D ( \AAZ)$ vanishes.
\edfn
\bppsn
A connection is torsionless iff its values on the generators $\sigma_1 ,\sigma_2 ,
\sigma_3 $ are  given by $$\nabla (\sigma_1)= \left( \matrix { \Box & a & b \cr a & \Box & c \cr
b & c & \Box \cr } \right),
\nabla (\sigma_2)= \left( \matrix { \Box & d & e \cr d & \Box & f  \cr
e & f & \Box \cr } \right) ,
\nabla (\sigma_3)= \left( \matrix { \Box & p-1  & q \cr  p & \Box &  r \cr
q & r & \Box \cr } \right).$$
\eppsn
{\it Proof:}
\bean
 \delta (\sum_{i,j} a_i \delta_j(b_i) \sigma_j)  & = &
            - \sqrt{-1} (  \sum_i ( \delta_1(a_i)   \delta_2(b_i)  -   \delta_2(a_i)     \delta_1(b_i) ) \sigma_3  \cr
& & +   \sum_i ( \delta_2(a_i)   \delta_3(b_i)  -   \delta_3(a_i)     \delta_2(b_i) ) \sigma_1  \cr
& & +  \sum_i ( \delta_3(a_i)   \delta_1(b_i)  -   \delta_1(a_i)     \delta_3(b_i) ) \sigma_2   )  \cr
\eean
\bean
m \circ \nabla (\sum_{i,j} a_i \delta_j (b_i) \sigma_j)
& =  &  m ( \sum_{i,j} \delta(   a_i \delta_j (b_i)) \otimes  \sigma_j )   +  \sum_{i,j}   a_i \delta_j (b_i) m \circ \nabla (\sigma_j )  \cr
   &  =  &    m ( \sum_{i,j,k} \delta_k (   a_i \delta_j (b_i)) \sigma_k \otimes  \sigma_j )  +  \sum_{i,j}   a_i \delta_j (b_i) m \circ \nabla  (\sigma_j)   \cr
\eean
Torsion of $\nabla$ vanishes iff $ (\delta - m \circ \nabla )( \sum a_i \delta_j (b_i) \sigma_j) \equiv 0$, or equivalently ,
\bean
 \sum_i ( \delta_j(a_i)   \delta_k(b_i)  -   \delta_k(a_i)     \delta_j(b_i) )   & =  &  \sum_i
 (\delta_j (   a_i \delta_k (b_i)) - \delta_k (   a_i \delta_j (b_i)))   \cr
   &  &+ \sum_{i,l}  a_i \delta_l (b_i) {(m \circ \nabla(\sigma_l))}_n \cr
\eean
whenever $ j \ne k$ and $n$ satisfies $ \sigma_j  \sigma_k  \sigma_n =\sqrt{-1}$
This happens iff
$$ 0= \sum_i a_i [ \delta_j, \delta_k](b_i) +    \sum_{i,l}  a_i \delta_l (b_i) {(m \circ \nabla(\sigma_l))}_n  $$
Using the Lie algebra relations between the $\delta_j$'s we get equivalence of the above system of equations with
\bean
0 & = &\sum_i a_i \delta_3 (b_i) +  \sum_{i,l}  a_i \delta_l (b_i) {(m \circ \nabla(\sigma_l))}_3   \cr
0 & = &    \sum_{i,l}  a_i \delta_l (b_i) {(m \circ \nabla(\sigma_l))}_2  \cr
0 & = &   \sum_{i,l}  a_i \delta_l (b_i) {(m \circ \nabla(\sigma_l))}_1 \cr
\eean
whenever $ j \ne k$ and $n$ satisfies $ \sigma_j  \sigma_k  \sigma_n =\sqrt{-1}$.\\
Taking $ b_i=\phi_{0,1},a_i=1$ we get  $ \delta_1(b_i)=  \delta_3(b_i)  =0,  \delta_2(b_i)=b_i$.
Putting these values in the above relations we get   ${(m \circ \nabla(\sigma_2))}_j=0$ for $ j=1,2,3$\\
Similarly taking $ b_i= \phi_{1,0},a_i = 1$ we get   $  {(m \circ \nabla(\sigma_1))}_j=0$ for $ j=1,2,3$\\
Substituting these values in the above equations we get,   $$  {(m \circ \nabla(\sigma_3))}_1=
{(m \circ \nabla(\sigma_3))}_2  =  \sum_i a_i \delta_3(b_i) ( 1+ (m \circ \nabla (\sigma_3))_3=0$$
Note,$ J= \{ \sum a_i \delta_3 ( b_i) | n \in \IN, a_1, \ldots ,a_n, b_1, \ldots , b_n \in \AAZ \}$
is a nontrivial ideal in $\AAZ$ hence equals $\AAZ$. Therefore $ {(m \circ \nabla (\sigma_3))}_3 = - 1 $.  Now the result follows from the anticommutation relation between the spin matrices.   \qed


\bdfn
A connection on a Finitely generated projective $\AAZ$ module ${\cal E}$, endowed with an $\AAZ$ valued inner product $ \langle  \cdot ,  \cdot \rangle$ is called unitary if$$
\delta \langle s,t \rangle = \langle \nabla  s,   t  \rangle  - \langle  s  , \nabla  t  \rangle,\forall s,t \in {\cal E} $$
Where the right hand side of this equation is defined by \\ $\langle \omega \otimes s , t \rangle= \omega \langle s, t \rangle, \langle s,\eta \otimes t \rangle = \langle s , t \rangle \eta^*$
\edfn
\bppsn
A connection $\nabla$ on ${\widetilde {\Omega}}^{1}_D ( \AAZ)$ is unitary iff   its values  on the generators
$\sigma_1 ,\sigma_2 ,
\sigma_3 $ are  given by $$\nabla (\sigma_1)= \left( \matrix { x & y & z \cr y  & u  & p \cr
z & v & q  \cr } \right),
\nabla (\sigma_2)= \left( \matrix { y  & u  & v  \cr u   & r  & s  \cr
p  & s  &  f  \cr } \right),
\nabla (\sigma_3)= \left( \matrix { z & p  & q \cr  v & s   &  f \cr
q & f  & g \cr } \right)$$
where $x,y,z,p,q,r,s,u,v,f,g \in \AAZ $ are selfadjoint elements.
\eppsn
{\it Proof:} Taking $s= a_i \sigma_i ,t= b_j \sigma_j$ in the defining condition of a  unitary connection we get
\bea \delta(\delta_{ij} a_i b_j^*) & = & a_i ( \langle \nabla( \sigma_i ), \sigma_j \rangle - \langle, \sigma_i , \nabla ( \sigma_j) \rangle) b_j^* \cr
& & + \delta_{ij} ( \delta(a_i) b_j^* - a_i {( \delta ( b_j))}^* ) \cr
\eea
implying $ \langle \nabla( \sigma_i ), \sigma_j \rangle = \langle  \sigma_i , \nabla ( \sigma_j) \rangle$\\
That is jth row of $\nabla(\sigma_i)$ is the star of the ith column of $\nabla(\sigma_j)$.
This completes the proof. \qed


\bcrlre
A connection $ \nabla$ can not simultaneously be torsionless and unitary.
\ecrlre
{\it Proof:} If possible let $\nabla$ be one such. Let $ v,p $ be as in proposition 3.16 and $c$ be as in proposition 3.14. Then $ v=c=p $   and also  $ v - p = -1 $.  This leads to a contradiction.   \qed


\newsection { connections with positive scalar curvature}
\bdfn
({\bf Theorem 2.9 of \cite {FGR} }) There is a sesquilinear map ${\langle \cdot,\cdot \rangle }_D : {\widetilde {\Omega}}^{k}_D ( \AAZ)   \otimes {\widetilde {\Omega}}^{k}_D ( \AAZ) \raro  {\bar {\cla_\hbar}}^w $
satisfying $ ( x, {\langle \omega, \eta \rangle }_D )= \int x \eta {\omega}^* $, \\  forall $x \in \cla_\hbar$
\edfn
In the following proposition we identify $ {\widetilde {\Omega}}^{k} ( \AAZ)$ with $ \AAZ \otimes M_2(\IC  )$
\bppsn
$ {\langle \omega , \eta \rangle}_D =\frac {1}{2} ( I \otimes tr) ( \omega  \eta^*)$
\eppsn
{\it Proof:}  Let $ \omega = \omega_0 \otimes I_2 + \sum_{i=1}^3 \omega_i \otimes \sigma_i$,
$ \eta = \eta_0 \otimes I_2 + \sum_{i=1}^3 \eta_i \otimes \sigma_i$\\
Then $ \frac {1}{2}( I \otimes tr) ({\omega}\eta ^*)= \sum_{i=0}^3{\omega}_i \eta_i ^*$\\
$ ( x, \sum_{i=0}^3  {\omega}_i\eta_i^*)= \tau (x  \eta_i {\omega}^*_i ) =(x,{\langle \omega , \eta \rangle}_D )$  forall $x \in \cla_\hbar$. This completes the proof since $ \cla_\hbar$ is dense in ${\widetilde \clh}^0$    \qed.\\

Notation :-- Let $ \omega \in {\Omega}^{1}_D ( \AAZ) $. Since $ K + \delta K$ is an ideal in
${\Omega}^{\bullet}_D ( \AAZ)   $  \\ $m :{\widetilde {\Omega}}^\bullet ( \AAZ) \otimes_{\AAZ}  {\widetilde {\Omega}}^\bullet ( \AAZ)
\raro   {\widetilde {\Omega}}^\bullet ( \AAZ) $  induces two maps denoted by the same symbol
$m :{\widetilde {\Omega}}^1_D ( \AAZ) \otimes_{\AAZ}  {\widetilde {\Omega}}^\bullet ( \AAZ)
\raro   {\widetilde {\Omega}}^{\bullet  +1} ( \AAZ) $\\
$m :{\widetilde {\Omega}}^\bullet ( \AAZ) \otimes_{\AAZ}  {\widetilde {\Omega}}^1_D ( \AAZ)
\raro   {\widetilde {\Omega}}^{\bullet+1} ( \AAZ) $. These inturn induce bounded maps $ m_L(\omega),
 m_R(\omega)  :{\widetilde \clh}^k \raro   {\widetilde \clh}^{k+1}$ as in remark (\ref {52}). \\
Since $ {\widetilde {\Omega}}^1_D ( \AAZ) $ is free, curvature of a connection $\nabla,$  \\  $ R ( \nabla)=-\nabla^2: {\widetilde {\Omega}}^1_D ( \AAZ) \raro {\widetilde {\Omega}}^2_D ( \AAZ)
\otimes_{\AAZ} {\widetilde {\Omega}}^1_D ( \AAZ) $ is given by a $3 \times 3$ matrix $ ((R_{ij}))$
with entries in ${\widetilde {\Omega}}^2_D ( \AAZ)$. Let $ P_{\delta K_1} :{\widetilde \clh}^2
\raro {\widetilde \clh}^1$ be the projection onto closure of $ \pi (\delta K_1) \subseteq    {\widetilde {\Omega}}^2_D ( \AAZ)$, and $ {R_{ij}}^\perp = (I- P_{\delta K_1})(R_{ij})$.Let $ e_1,e_2,e_3$ be the canonical basis of ${\widetilde {\Omega}}^1_D ( \AAZ)   $.  If we denote by $ {Ric}_j= \sum_i { m_L(e_i)}^{ad} (  {R_{ij}}^\perp)   \in {\widetilde \clh}^{1}  $  then Ricci curvature of $\nabla$ is given by
$$ Ric(\nabla) = \sum_{j} Ric_j \otimes e_j \in {\widetilde \clh}^{1} \otimes_{\AAZ}  {\widetilde {\Omega}}^1_D ( \AAZ)$$
Here superscript ad stands for Hilbert space adjoint. Finally the scalar curvature $ r ( \nabla)$ of $\nabla$ is given by  $$ r(\nabla )= \sum_i {m_R(e_i^*)}^{ad}  (Ric_i)  \in  {\widetilde \clh}^{0} $$
\bppsn
Let $f,g :\IT \raro \IR $ be smooth maps. Henceforth we visualize $f,g $ as elements of $ S^c$ in the following way, $ f(x,y,p)=\delta_{0p} f(x),$  \\ $ g(x,y,p)= \delta_{0p} g(y)$. Similar considerations will be applied for $f^\prime,g^\prime$.
Let $\nabla$ be the connection given by $ \nabla ( \sigma_1)=f^\prime \delta(g) \sigma_1 +g^\prime \delta(f) \sigma_2, \nabla ( \sigma_2)=g^\prime \delta(f) \sigma_1, \nabla ( \sigma_3)=0$, then $r(\nabla)$ is $ -2
f^{\prime 2},g^{\prime 2} $
\eppsn
{\it Proof:}   By direct computation one gets,
$$\nabla^2(\sigma_1)= -R_{11} \sigma_1 - R_{12} \sigma_2 ,  \nabla^2(\sigma_2)=-R_{21} \sigma_1,
\nabla^2(\sigma_3)=0$$ where $$ R_{11}= f^{\prime \prime} g \sigma_3 ,  R_{12}= \sqrt {-1} (f^{\prime 2} g^{\prime 2} -g^{\prime \prime}f^\prime) \sigma_3,
 R_{21} =- \sqrt {-1} (g^{\prime \prime}f^\prime  +f^{\prime 2} g^{\prime 2}  )\sigma_3,$$ other $ R_{ij}$'s are zero. \\Then $$ Ric_1= -f^{\prime \prime} g \sigma_2 -(
g^{\prime \prime}f^\prime  +f^{\prime 2} g^{\prime 2}  )\sigma_1,
 Ric_2= (g^{\prime \prime}f^\prime  -f^{\prime 2} g^{\prime 2}  )\sigma_2  $$
implying the desired conclusion
$ r(\nabla)=  -2  f^{\prime 2}g^{\prime 2} $.  \qed
\brmrk
(i) All these notions of Ricci curvature, scalar curvature was introduced by \cite {FGR}. To the best of our knowledge it is the first infinite dimensional example where one can have connections with nontrivial scalar curvature.\\
(ii) Note out choice of the spectral triple depend on a parameter $\al$. For the connections we have considered scalar curvature does not depend on the parameter $\al$.
\ermrk
\newsection {nontriviality of the chern character associated with the spectral triples}
The spectral triple we constructed depends on a real parameter $\al$. In this section we show that the Kasparov module associated with the spectral triple \cite {CON} \cite {BLA}  are homotopic. We also argue that they give non-trivial elements in $K^1( \cla_\hbar)$ by explicitly computing pairing with some unitary in the algebra representing elements of $ K_1(\cla_\hbar)$.
\blmma
Let A be a selfadjoint operator with a bounded inverse and B a symmetric operator with $ D(A) \subseteq D(B)$ on some Hilbert space $\clh $ . Also suppose that $\| B u \| \le  a \| A u \| , \forall u \in D(A) $. Then $ {| A|}^{-p} B {|A |}^{-(1-p)} \in  \clb ( \clh )$ and $ \|  {| A|}^{-p} B {|A |}^{-(1-p)}\| \le a $.
\elmma
{\it Proof :} Clearly $ \| B u \| \le a \| |A| u \| , \forall u \in D(A)$  implying $ \| B {| A | }^{-1} \| \le a $.  For $ u,v \in D(A)$
\bean
\| {|A|}^{-1} B u \| & = & \sup_ {v \in D(A), \| v \| \le 1} | \langle {| A | }^{-1} B u , v \rangle | \cr
  & = & \sup_ {v \in D(A), \| v \| \le 1}  | \langle u,    B   {| A | }^{-1}  v   \rangle |  \le a \| u \| \cr
\eean
Therefore   $ {| A | }^{-1} B \in \clb ( \clh ),  \| {| A | }^{-1} B\| \le a $ \\
Let $ \clh_p $ be the Hilbert space completion of $ \cap D ( A^n) $ with respect to \\  $ {\| u \| }_p= \| {| A |}^p u \|.$ Let $ B_1 : \clh_1 \raro \clh_0, B_0 : \clh_0 \raro \clh_{-1}$ be the maps given by $ B_i (u)= B(u) $ for $ u \in \cap D ( A^n)$. Then $  \| B_ 1 \|  ,   \| B_ 1 \|   \le   a $.  By Calderon-Zygmund interpolation theorem \cite {RS}  we get maps $ B_p : \clh_p \raro \clh_{-(1-p)}$ for $ 0 \le p \le 1$ with $ \| B_p \| \le a $. On  $ \cap D ( A^n), B_p $ agrees with  ${| A|}^{-p} B {|A |}^{-(1-p)}$  proving the lemma.   \qed  \\
\blmma
\label {72}
Let $A,B$ be as above with $ a < 1$. Let $ A_t = A + t B, t \in [0,1]$. Then $ t \mapsto {\tan}^{-1} (A_t) $ is a norm continuous function.
\elmma
{\it Proof:} Let $ C=  {| A|}^{-1/2} B {|A |}^{-1/2}  $ , then by the previous lemma $ \| C \| \le a $.
For $ \lmd \in i \IR , \| | A | {( A - \lmd )}^{-1} | \le 1$.
\bean
A_t - \lmd & = &  ( A - \lmd ) + t {|A|}^{1/2} C  {|A|}^{1/2} \cr
  & = & {|A|}^{1/2}( ( A - \lmd ) {| A|}^{-1}+ t C  ){|A|}^{1/2}  \cr
  & = &  {| A|}^{1/2} ( 1 + t C  {( A - \lmd )}^{-1} | A| )  ( A - \lmd ) {| A|}^{-1}{| A|}^{1/2}  \cr
\eean
Now note $ \|  t C  {( A - \lmd )}^{-1} | A|  \| \le  a < 1$ for $ 0 \le t \le 1 $. Therefore \\
${(A_t - \lmd )}^{-1} = {|A |}^{-1/2} | A|   {( A - \lmd )}^{-1}{ ( 1 + t C  {( A - \lmd )}^{-1} | A|  )}^{-1} {|A |}^{-1/2}$ \\
So, if we denote by $ R_t (\lmd)={(A_t - \lmd )}^{-1}$ and $ F( \lmd)= | A|   {( A - \lmd )}^{-1}$ then the above equality becomes,
\bea
\label {14}
 R_t (\lmd) & = &  {|A |}^{-1/2} | A|  R_0  (\lmd) {( 1 + t C | A| F ( \lmd) )}^{-1}{|A |}^{-1/2} \cr
  & =  & R_0  (\lmd) +  {|A |}^{-1/2}  F ( \lmd)  \sum_{n=1}^\infty  {(-tC  F ( \lmd))}^n  {|A |}^{-1/2}  \cr
\eea
Let $ \lmd \in \IR ,t,s \in [0,1],u,v \in D(A) $. Observe

(i) $$ \|  \sum_{n=1}^\infty  {(-tC  F ( i \lmd))}^n  {|A |}^{-1/2}  u -\sum_{n=1}^\infty  {(-sC  F ( i  \lmd))}^n  {|A |}^{-1/2}  u \| $$
$$ \le   \sum_{n=0}^\infty  \|  (t^{n+1} -s^{n+1}) { C  F (  i \lmd)  \|}^n  \|C \|  \|F (  i \lmd)  {|A |}^{-1/2}  u\|  $$
$$ \le  \sum_{n=0}^\infty |(t^{n+1} -s^{n+1}) | a^n a \|F (  i \lmd)  {|A |}^{-1/2}  u\| $$
$$ \le |(t-s)| \sum_{n=0}^\infty  n+1 a^{n+1}   \|F (  i \lmd)  {|A |}^{-1/2}  u\|   $$    $$ \le |(t-s)| \frac {a} {{(1-a)}^2} \|F (  i \lmd)  {|A |}^{-1/2}  u\|  $$

(ii)
\bean
\int_0^\infty  {\|F (  i \lmd)  {|A |}^{-1/2}  u\| }^2 d \lmd &  \le &   \int_0^\infty \langle {(A^2+ \lmd^2)}^{-1} u,  | A | u \rangle d \lmd  \cr
& = & \frac {1}{2} \int_0^\infty \langle {(A^2+ \xi )}^{-1} u,  | A | u \rangle \frac {d  \xi} { \sqrt \xi} \cr
  & = & \frac {1} {2} \pi \langle {A^2}^{-1/2 } u,  | A | u \rangle = \frac {\pi }{2} { \| u \| }^2  \cr
\eean

(iii) Using (\ref {14}), (i) ,(ii)  we get
$$    \int_0^\infty |\langle (R_t( i \lmd ) -R_s( i \lmd )) u, v \rangle | d \lmd $$
$$\le  \int_0^\infty  |(t-s)| \frac {a} {{(1-a)}^2} \|F (  i \lmd)  {|A |}^{-1/2}  u\|
\|F ( - i \lmd)  {|A |}^{-1/2}  v\| d \lmd  $$
$$ \le   |(t-s)| \frac {a} {{(1-a)}^2}{(\int_0^\infty  {\|F (  i \lmd)  {|A |}^{-1/2}  u\| }^2 d \lmd )}^{1/2}
{(\int_0^\infty  {\|F (-  i \lmd)  {|A |}^{-1/2}  v \| }^2 d \lmd )}^{1/2} $$
$$\le |(t-s)| \frac {a} {{(1-a)}^2} \frac{\pi}{2} \|u \| \| v \| $$
This shows $ \lim_{s \raro t} \| \int_0^\infty  (R_t( i \lmd ) -R_s( i \lmd ))  d \lmd  \| = 0 $. Similarly one can show   $ \lim_{s \raro t} \| \int_0^\infty  (R_t( -i \lmd ) -R_s(-i \lmd ))  d \lmd  \| = 0 $.  Now the result follows once we observe $ {\tan}^{-1} A_t = \int_0^\infty   ( R_t( i \lmd ) + R_t ( -i \lmd)) d \lmd $.  \qed\\
\blmma
\label {73}
Let $A,B$ be as above except now we do not require $A$ to be invertible. Instead we assume $A$ to have discrete spectrum. Then there exists $ \kappa  \ge 0$ such that $ t \mapsto {\tan}^{-1} ( A_t + \kappa )$ is norm continuous.
\elmma
{\it Proof:}  Without loss of generality we can assume 0 is an eigenvalue of $A$. Otherwise we are done by the previous lemma.
Choose $2 \le n \in \IN $ such that  $ b  =a \frac {n} { n-1} < 1$  Choose $\kappa > 0 $ such that  \\ (i) smallest positive eigenvalue of $A$ is greater than $\kappa$.  \\ (ii) if $\beta$ is the biggest negative eigenvalue then $ \beta < n \kappa$.\\Let $ {\widetilde  A} = A+ \kappa, {\widetilde A_t}= {\widetilde  A}  +t B $. Then by choice of $\kappa$\\ (i) ${\widetilde  A}  $  is an invertible selfadjoint operator.  \\
(ii) $ \| B { {\widetilde  A}}^{-1} \| \le  a \| A {(    A+ \kappa)}^{-1} \| \le   a \frac {n} { n-1}  < 1$  \\
That is $B$ is relatively bounded with respect to $   {\widetilde  A}  $  with relative bound $b < 1$.
Now an application of the previous result to the pair ${\widetilde  A},B$  does the job.  \qed\\
Combining these two we get
\bppsn
Let $A,B$ be operators on the Hilbert space $ \clh$ such that \\
(i) A is selfadjoint with compact resolvent.\\(ii) B is symmetric with $ D(A) \subseteq D( B)$, and relatively bounded with respect to A with relative bound less than 1. \\
Then there exists  a continuous function $f : \IR \raro \IR $ satisfying \\$ \lim_{x \raro \infty} f(x)=1,
\lim_{x \raro -\infty} f(x)=-1$ such that $ t \mapsto f( A + t B ) $ is norm continuous.
\eppsn
{\it Proof:} If $ A$ is invertible then by lemma (\ref {72}) $f(x) = \frac {2} {\pi} {\tan}^{-1} (x) $ serves the purpose. In the other case by lemma (\ref {73}) $ f(x) =\frac {2} {\pi} {\tan}^{-1} (x  +  \kappa) $  does the job. \qed\\
Let the Hilbert space $ \clh$  and the operators $A,B,D$ be as in corollary (\ref {33}.
\bcrlre
\label {EQ}
The Kasparov module associated with $ ( \AAZ ,\clh, D)$ is operatorial homotopic with $ ( \AAZ ,\clh, A)$
\ecrlre
{\it Proof:}  Let $A_t=A+t B$ for $t \in [0,1]$. Then $D=A_1,A=A_0$.  As remerked earlier $  ( \AAZ ,\clh, A_t)$   are spectral triples. Let $f$ be the function obtained from the previous proposition for the pair $A,B$. Then ${((\cla_\hbar, \clh,   f(A_t)))}_{t \in [0,1]} $ gives the desired homotopy. \qed.

As remarked earlier the operator $A$ depends on a real parameter $\al >1$. Now we will make that explicit and denote $A$ by ${A}^{(\al)}$.
\bppsn
\label {76}
The Kasparov modules associated with $ ( \AAZ ,\clh,{A}^{(\al)})$  are operatorially homotopic for $ \al > 1$
\eppsn
{\it Proof:}  By proposition (\ref {26}), $ \clh = L^2( \IT \times \IT \times \IZ) \otimes {\IC}^2$. Let $B$ be the operator $-2\pi c M_p \otimes \sigma_3$. Here p denotes the $\IZ$ variable in the $L^2$ space. Then $B$ is selfadjoint with $ D({A}^{(\al)}) \subseteq D(B)$. Also $B$ is relatively bounded with respect to ${A}^{(\al)}$ with relative bound less than $\frac {1} { \al}  < 1  $.Let $ {  {A}^{(\al)}  }_t=
{A}^{(\al)}  + t B $for $ t \in [0,1]$. Then   $  {  {A}^{(\al)}  }_t  =    {A}^{(\al+t)}$.  Let $f$ be the function obtained from proposition 5.4 for the pair ${A}^{(\al)},B$. Then from the norm continuity of $ t \mapsto f({A}^{(\al+t)})  $ we see the Kasparov modules   $ {(( \AAZ ,\clh,  {A}^{(\al+t)}))}_{t \in [0,1]}
$ are homotopic. Since $\al$ is arbitrary this completes the proof.  \qed
\brmrk
Proposition ( \ref {76}) and corollary (\ref {EQ}) together imply the Kasparov module associates with the spectral triple
$( \AAZ ,\clh, D)$ is independent of $ \al $.
\ermrk
In the next proposition we show $ (( \AAZ ,\clh, D)$ has non trivial chern character.
\bppsn
The Kasparov module associated with $(( \AAZ ,\clh, D)$ gives a nontrivial element in $ K^1 ( \cla_\hbar)$
\eppsn
{\it Proof:} By corollary \ref {EQ} $(( \AAZ ,\clh, D)$ and $ ( \AAZ ,\clh, A)$ give rise to same element
  $ [(\AAZ, \clh, A )] \in   K^1 ( \cla_\hbar)$.  Let $ \phi \in \AAZ$ be the unitary whose symbol in $S^c$ is given by $ \phi (x,y,p) = \delta_{0p} e^{2\pi i y}$. This gives an element $[ \phi ] \in K_1(\cla_\hbar )$.
It suffices to show $ \langle  [ \phi ]  ,   [(\AAZ, \clh, A )] \rangle \ne 0 $ where $   \langle  \cdot ,  \cdot \rangle  : K_1( \cla_\hbar) \times K^1( \cla_\hbar)  \raro \IZ $ denotes the pairing   coming from Kasparov product. $\phi $ acts on $ L^2 (\cla_\hbar ) \otimes \IC^2  \cong L^2( [0,1] \times \IT \times \IZ)  \otimes \IC^2$ as a composition of two  commuting unitaries unitaries $U_1 = M_{e(y)} \otimes I_2, U_2= M_{e(p \nu \hbar )}\otimes I_2$. Then note $ U_2 $ commutes with $ A$. Let $E$ be the projection $ E=I ( A \ge 0)$.  $U_2$ also commutes with $E$. Now  by proposition 2 ( page 289 of \cite {CON})  $ E U_1 U_2 E$ is a Fredholm operator and  $ \langle  [ \phi ]  ,   [(\AAZ, \clh, A )] \rangle = Index( E U_1 U_2 E)  =Index( E U_1  E)$, last equality holds because $U_2$ commutes with $E$. Now $Index( E U_1  E) \ne 0$ because this is the index pairing of the Dirac operator on $ \IT^3$ with the unitary $U_1$.  \qed.

  \end{document}